\documentclass[notitlepage,11pt,reqno]{amsart}
\usepackage[hmargin={1.2in, 1.2in}, vmargin={1.3in, 1.3in}]{geometry}
\usepackage{amssymb,amsmath,amsthm}
\usepackage[mathscr]{eucal}
\usepackage[breaklinks=true]{hyperref}

\newcommand{\F}{\ensuremath{\mathbb{F} }}

\newcommand{\Z}{\ensuremath{\mathbb{Z} }}

\newcommand{\m}{\ensuremath{\mathfrak{m} }}

\newcommand{\tangle}[1]{\ensuremath{\langle #1 \rangle}}
\newcommand{\commen}[1]{}

\DeclareMathOperator{\Ann}{\ensuremath{Ann}}

\theoremstyle{plain}
\newtheorem{theorem}[equation]{Theorem}

\newtheorem{prop}[equation]{Proposition}
\newtheorem{cor}[equation]{Corollary}

\theoremstyle{definition}
\newtheorem{definition}[equation]{Definition}
\newtheorem{remark}[equation]{Remark}

\numberwithin{equation}{section}

\makeatletter
\@namedef{subjclassname@2020}{%
  \textup{2020} Mathematics Subject Classification}
\makeatother

\begin{document}
\title{Covering modules by proper submodules}

\author{Apoorva Khare}
\address[A.~Khare]{Indian Institute of Science, Bangalore 560012, India;
and Analysis and Probability Research Group; Bangalore 560012, India}
\email{\tt khare@iisc.ac.in}

\author{Akaki Tikaradze}
\address[A.~Tikaradze]{Department of Mathematics, University of Toledo,
Toledo 43606, USA}
\email{\tt tikar06@gmail.com}


\date{\today}

\begin{abstract}
A classical problem in the literature seeks the minimal number of proper
subgroups whose union is a given finite group. A different question, with
applications to error-correcting codes and graph colorings, involves
covering vector spaces over finite fields by (minimally many) proper
subspaces. In this note we cover $R$-modules by proper submodules for
commutative rings $R$, thereby subsuming and recovering both cases above.
Specifically, we study the smallest cardinal number $\aleph$, possibly
infinite, such that a given $R$-module is a union of $\aleph$-many proper
submodules.
(1)~We completely characterize when $\aleph$ is a finite cardinal; this
parallels for modules a 1954~result of Neumann.
(2)~We also compute the covering (cardinal) numbers of finitely generated
modules over quasi-local rings and PIDs, recovering past results for
vector spaces and abelian groups respectively.
(3)~As a variant, we compute the covering number of an arbitrary direct
sum of cyclic monoids.
Our proofs are self-contained.
\end{abstract}
\maketitle

\section{Introduction}

\noindent The covering problem is a well-known and classical question in
group theory:\smallskip

\textit{Given a non-cyclic group $G$, find the minimal (cardinal) number
of proper subgroups, termed the \textbf{covering number} $\sigma(G)$,
whose union is $G$.}\smallskip

This problem has a long history in the literature. It is immediate that
two proper subgroups can never cover $G$; in 1926,
Scorza~\cite{Scorza} showed that $\sigma(G) = 3$ if and only if $G$
has a quotient isomorphic to the Klein-4 group $(\Z / 2 \Z)^2$. 
This question and related variants have since been the subject of a vast
number of papers; a small sampling from the 20th century is
\cite{Br2,Br3,Cohn,HR,Neu,RaoReid,Rosenfeld,Tom1,Tom2,Tom3}.
The subject also continues to attract much attention more recently, see
e.g.\ \cite{AS,BLR,Ber,Bh1,Bh2,Bh3,BFM,Br1,FS,LG,Mar,Sun,Swartz,Sz} and the
references therein.

One of the early results on (finite) coverings of groups is by Neumann
(1954), who worked in slightly greater generality -- for unions of
\textit{cosets} of subgroups:

\begin{theorem}[Neumann, \cite{Neu}]\label{Tneumann}
Suppose a group $G$ is the union of finitely many cosets
\[
g_1 G_1, \dots, g_n G_n, \quad \text{where } g_i \in G, \ G_i \leq G,
\]
where we write $H \leq G$ to denote that $H$ is a subgroup of $G$.
\begin{enumerate}
\item Then the index $[G : G_i]$ is finite for some $i$.

\item Now suppose no proper sub-union equals $G$. Then $[G : G_i]$ is
finite for all $i$. In particular, if moreover $g_i = e\ \forall i$, then
$[G : G_1 \cap \cdots \cap G_n]$ is finite.

\item Consequently, $G$ can be written as a union of (i.e., admits a
covering by) finitely many proper subgroups, if and only if $G$ has a
finite non-cyclic quotient.
\end{enumerate}
\end{theorem}

On a tangential note, we remind the reader that covering the specific
group $G = \Z$ by finitely many cosets (with distinct moduli) was a
subject well-studied by Erd\H{o}s; for instance, some of his conjectures
(including with Graham) were settled in~2007 -- in a stronger form -- by
Filaseta--Ford--Konyagin--Pomerance--Yu~\cite{FFKPY}.

Returning to coverings of general groups, Rao--Reid used Neumann's
theorem~\ref{Tneumann}, together with ideas of Rao--Rao~\cite{RR}, to
prove:

\begin{theorem}[Rao--Reid \cite{RaoReid}]\label{Traoreid}
Suppose an abelian group $G$ admits an irredundant covering by finitely
many (proper) subgroups.
\begin{enumerate}
\item Every such covering is induced by lifting a covering of a finite
quotient of $G$.

\item The covering number $\sigma(G) = p+1$, where $p \geq 2$ is the
smallest prime such that $G/pG$ is not cyclic.
\end{enumerate}
\end{theorem}

It is natural to try and extend such results to other settings. For
instance, coverings of rings have been studied (including very recently)
in~\cite{BBS,BBBK,CaiWer,LM,PW,Werner}.
Yet another motivation comes from the short note~\cite{Khvec}, where the
first named author found a sharp bound for the number of proper subspaces
of a fixed codimension $d \geq 1$ needed to cover a vector space. We
present here the $d=1$ result -- which is folklore, and is all that is
required for the purposes of this paper:

\begin{prop}\label{Pvecunion}
Suppose $V$ is a vector space of dimension at least $2$ over a field
$\F$, and define $\nu(\F,V)$ to be the cardinal number
\begin{equation}\label{Enu1}
\nu(\F,V) := \begin{cases}
|\Z|, \qquad & \text{if } \F,
\dim V \text{ are both infinite,}\\
|\F|+1, & \text{otherwise.}
\end{cases}
\end{equation}
Now $V$ can be written as a union of $\aleph$-many proper vector
subspaces, if and only if $\aleph \geq \nu(\F,V)$ as cardinal numbers.
In other words, $\sigma(V) = \nu(\F,V)$.
\end{prop}

(Here, $\sigma(V)$ of course denotes the covering number of $V$ in the
category of vector spaces.) As we show in Corollary~\ref{Clocal} below,
the above result extends to cover free modules $R^n, n \geq 2$ over a
local ring $(R,\m)$: we prove $\sigma(R^n) = |R/\m|+1$ as cardinal
numbers, irrespective of whether or not $R/\m$ is finite.
For now, we provide a quick proof of Proposition~\ref{Pvecunion} when
$\F$ is finite. That $\aleph \geq |\F|+1$ works, follows by lifting the
cover by lines of the plane $\F^2$, i.e., projective space, to $V$.
Conversely, if $\{ V_i : i \in \mathscr{I} \}$ is a finite and
minimal/irredundant cover of $V$, we may suppose each $V_i$ has
codimension one, whence $\cap_i V_i$ has finite codimension in $V$.
Working modulo $\cap_i V_i$ now reduces the case to $\dim V < \infty$, so
$|V| < \infty$ and $|V_i| = |V|/|\F|$. It follows that $|\mathscr{I}| =
\sigma(V) > |\F|$.\medskip

Covering vector spaces by proper subspaces, and variants of this problem,
find applications to error-correcting codes and block designs over finite
fields (see the references in \textit{loc.~cit.}), as well as to graph
colorings (see e.g.~\cite{Sz}). Other variants involving covering a
vector space by (translates/shifts) of subspaces can be found
in~\cite{CER,Jam}.

\subsection{Contributions}

The goal of this note is to unify and extend the study of coverings of
abelian groups and of vector spaces by working with $R$-modules, where
$R$ will always denote a unital commutative ring. We list here a few of
the contributions.
\begin{enumerate}
\item We completely characterize the modules over an arbitrary ring $R$,
that admit a finite covering. For $R = \Z$, we recover the result for
finite coverings of abelian groups; our proof is self-contained and works
for arbitrary $R$.
\item We next determine the covering number in a large number of cases,
including: (a)~all finitely generated modules over quasi-local rings;
(b)~all divisible modules over a PID;
(c)~direct sums of cyclic modules over a local ring (which subsumes the
case of vector spaces).
\end{enumerate}

\section{Finite coverings of modules, and coverings of finitely generated
modules}

In this section, we extend the results for abelian groups above, first to
finite coverings of (arbitrary) $R$-modules $M$ for arbitrary unital
commutative rings $R$. This is quickly followed by the case of coverings
of finitely generated modules over quasi-local rings (i.e., ones with
only finitely many maximal ideals). Our proofs are self-contained.

Our first result completely characterizes -- for an arbitrary unital
commutative ring $R$ -- when an (arbitrary) $R$-module $M$ admits a
finite covering; note this extends Theorem~\ref{Traoreid}. In particular,
our proof (of the nontrivial, second part) differs from that
in~\cite{RaoReid}, as this latter proof for groups does not extend to
$R$-modules.

We begin with two easy observations.
First, if $M$ is a finite set and an $R$-module, and $\m M \subsetneq M$
for a maximal ideal $\m \subset R$, then $R / \m$ is a finite field
(since $M / \m M$ is a finite-dimensional vector space over $R / \m$).
From this it follows that $\m$ contains a (unique) prime integer $p_\m
\in \Z$, i.e., $\m \cap R_0 = p_\m R_0$, where $R_0$ is the unital
subring of $R$ generated by $1_R$.
Second, from this it follows that
\begin{equation}\label{Eprimedec}
M = \oplus_{p \in \Z \text{ prime}} M_{(p)}
\end{equation}
as $R$-modules, where $M_{(p)}$ is the $p$-(torsion) subgroup of $M$ under
addition. We now have:

\begin{theorem}\label{Tfinite}
Suppose $R$ is a unital commutative ring and $M$ an $R$ module.
\begin{enumerate}
\item Let $S$ be the set of maximal ideals $\m$ for which $\dim_{R/\m}
M/\m M \geq 2$, and define
\begin{equation}\label{Eupperbound}
\nu'(R,M) := \min_{\m \in S} |R/\m|+1,
\end{equation}
as cardinal numbers. If $S$ is nonempty, then $\sigma(M) \leq \nu'(R,M)$.

\item $M$ admits a finite covering by proper submodules, if and only if
$M$ has a finite non-cyclic quotient $R$-module. Moreover, every such
irredundant finite covering is induced by lifting a covering of a finite
quotient of $M$.

\item If $M$ admits a finite covering by proper $R$-submodules, then the
covering number $\sigma(M) = \nu'(R,M)$.
\end{enumerate}
\end{theorem}

In other words, the definition of the cardinal number $\nu(\F,V) =
\nu'(\F,V)$ in Equation~\eqref{Enu1} can be extended from
\textit{finite-dimensional} $\F$-vector spaces $V$, to $R$-modules that
admit a finite covering: $\sigma(M) < \infty \implies \sigma(M) =
\nu'(R,M)$. Note, this cardinal equals $|R/\m_0|+1$ with $R/\m_0$ a
finite field.

To prove Theorem~\ref{Tfinite}, we first provide a short proof -- using
amenability -- of Neumann's theorem~\ref{Tneumann} for abelian groups.

\begin{prop}\label{Pneumann}
Suppose an amenable group $G$ is an irredundant union of finitely many
proper subgroups $G_1, \dots, G_n$. Then $[G : G_1 \cap \cdots \cap G_n]$
is finite (whence so is each $[G : G_i]$).

Consequently, an abelian group $G$ can be written as a union of (i.e.,
admits a covering by) finitely many proper subgroups, if and only if $G$
has a finite non-cyclic quotient.
\end{prop}

\begin{proof}
Since $G$ is amenable, denote the relevant finitely additive
left-invariant probability measure on $G$ by $\mu$.
By relabeling, there exists $0 \leq t \leq n$ such that $G_1, \dots, G_t$
have finite index in $G$, while $G_{t+1}, \dots, G_n$ do not. Note that
$\mu(G_i) = 0$ for $i>t$, so we must have $t>0$, whence the index
\[
m := [G : \cap_{i=1}^t G_i] \in (1,\infty).
\]
Now $G \setminus \cup_{i=1}^t G_i$ is a disjoint union of cosets of
$\cap_{i=1}^t G_i$, so if $t<n$, then we compute:
\[
1 - \frac{1}{m} = \mu \left( G \setminus \cup_{i=1}^t G_i \right) \leq
\mu(G_{t+1} \cup \cdots \cup G_n) \leq \sum_{i=t+1}^n \mu(G_i) = 0.
\]
This contradiction shows the first part. The second part is standard: if
$G$ is abelian and admits a finite covering, then so does its finite
quotient $G / \cap_{i=1}^n G_i$ from above, whence this cannot be a
cyclic group. The converse is immediate by lifting a cover by the cyclic
proper subgroups, of $G/H$ (a finite, non-cyclic quotient group).
\end{proof}

The next argument is useful in multiple situations, hence is isolated
into a standalone result. In it and in the sequel, $\sigma(M)$ will
always denote the covering number of an $R$-module $M$ by proper
$R$-submodules.

\begin{prop}\label{Pproduct}
Fix an integer $n \geq 1$, finitely many unital commutative rings $R_1,
\dots, R_n$, and an $R_k$-module $M^{(k)}$ for all $k$, such that
$M^{(k_0)}$ admits a covering by proper $R_{k_0}$-submodules for at least
one $k_0 \in [1,n]$. Letting $\mathscr{K}_0$ denote the set of all such
$k_0$, and setting $R := \oplus_{k=1}^n R_k,\ M := \oplus_{k=1}^n
M^{(k)}$, we have $\sigma(M) = \min_{k_0 \in \mathscr{K}_0}
\sigma(M^{(k_0)})$ as cardinal numbers.
\end{prop}

This minimum cardinal number exists (as do all such minima below)
e.g.~by~\cite{MeSa}.

\begin{proof}
That $\sigma(M) \leq \min_{k_0 \in \mathscr{K}_0} \sigma(M^{(k_0)})$ is
immediate: let $k' \in \mathscr{K}_0$ attain this minimum, and lift a
cover of $M^{(k')}$ to $M$.

Conversely, suppose $M = \cup_{i \in \mathscr{I}} M_i$ is a covering with
$|\mathscr{I}| = \sigma(M)$. We first get rid of the `cyclic' factors.
Namely, suppose without loss of generality that $M^{(k)}$ is cyclic as an
$R_k$-module, say $M^{(k)} = R_k v_k$, for $k=1, \dots, t$. Note that
$t<n$ since $\mathscr{K}_0$ is nonempty. Now given $M = \cup_{i \in
\mathscr{I}} M_i$, we claim that the module $M / (M^{(1)} \oplus \cdots
\oplus M^{(t)})$ over the ring $R' := \times_{k > t} R_k$ is covered by
$M_i / (M^{(1)} \oplus \cdots \oplus M^{(t)})$ for all $i \in
\mathscr{I}$ such that $v_1, \dots, v_t \in M_i$. (Call this subset
$\mathscr{I}_0 \subset \mathscr{I}$.)
Indeed, given $(v_k)_{k>t} \in M / (M^{(1)} \oplus \cdots \oplus
M^{(t)})$, we must have ${\bf v} := (v_1, \dots, v_t, (v_k)_{k>t}) \in
M_i$ for some $i \in \mathscr{I}$. But then $v_k = 1_{R_k} {\bf v} \in
M_i$ for $1 \leq k \leq t$, whence $i \in \mathscr{I}_0$; and
$\oplus_{k>t} 1_{R_k} \cdot {\bf v} = (v_k)_{k>t} \in M_i$ as well,
showing the claim.

Thus, we may assume now that every $M^{(k)}$ admits a covering by proper
$R_k$-submodules. Given a module $N$ over $R = \times_{k=1}^n R_k$, write
$N = \oplus_{k=1}^n N^{(k)}$, with $N^{(k)} = R_k N$. Thus for each $i$,
we have $M_i = \oplus_{k=1}^n M_i^{(k)}$ is proper, whence there exists
$1 \leq k_i \leq n$ with $M_i^{(k_i)} \neq M^{(k_i)}$. We may then
replace $M_i$ by
\[
M'_i := M_i^{(k_i)} \oplus \oplus_{k \neq k_i} M^{(k)},
\]
and still obtain a covering of $M$ by $\sigma(M)$-many proper
$R$-submodules $M'_i$.

Finally, let $\mathscr{I}_k := \{ i \in \mathscr{I} : k_i = k \}$ for $1
\leq k \leq n$. We claim there exists $1 \leq k \leq n$ such that $\{
M_i^{(k)} : i \in \mathscr{I}_k \}$ is a covering of $M^{(k)}$; notice
that this would imply that $\sigma(M) \geq \sigma(M^{(k)})$, and conclude
the proof. This claim is shown by contradiction: if no such $k$ exists,
then for every $k$ there exists $m_k \in M^{(k)}$ that does not belong to
$M_i^{(k)}$ for any $i \in \mathscr{I}_k$. But then $\oplus_{k=1}^n m_k
\not\in M'_i$ for all $i \in \mathscr{I}$.
\end{proof}

With these preliminary results at hand, we can now completely
characterize when an $R$-module $M$ has a finite covering.

\begin{proof}[Proof of Theorem~\ref{Tfinite}]\hfill
\begin{enumerate}
\item This is immediate: if $\dim_{R/\m} M/\m M \geq 2$, then by
Proposition~\ref{Pvecunion}, the quotient module (over the field $R /
\m$)
\[
M \twoheadrightarrow M / \m M \twoheadrightarrow (R / \m)^2
\]
is a union of $|R / \m| + 1$ proper submodules, which can then be lifted
to a covering of $M$. Now take the minimum over $\m \in S$.

\item First suppose $M$ admits a finite covering, say a
minimal/irredundant one, by proper submodules $M_1, \dots, M_n$. Then $M
/ \cap_i M_i$ is finite by Proposition~\ref{Pneumann}, and is not cyclic,
else the cyclic generator would be contained in some $M_{i_0} / \cap_i
M_i$. Moreover, the covering of $M$ by the $M_{i_0}$ is the lift of the
covering of $M / \cap_i M_i$ by the $M_{i_0} / \cap_i M_i$. Conversely,
if $M / N$ is finite and non-cyclic, then $M/N$ is the union of its
(finitely many) cyclic submodules, and this lifts to a finite covering of
$M$.

\item One inequality follows from~\eqref{Eupperbound}. To show the
reverse inequality, given a module $M$ admitting a finite covering, we
begin with two reductions.
(a)~From the preceding part, we may assume $M$ is finite.
(b)~We reduce to the case of finite rings, by working over $R_0 := R /
\Ann_R(M)$ instead of $R$. Indeed, $R_0$ embeds into ${\rm End}(M)$ (via
$r \mapsto (m \mapsto rm)$), which is a finite set.

Now we proceed. Given an irredundant covering $M = M_1 \cup \cdots \cup
M_n$, each $M_i$ has a nonzero simple quotient, say $R_0 / \m_i$, whence
$M_i / \m_i M_i \neq 0$. Let the ideal $I \subset R_0$ be the product of
the distinct ideals among $\m_1, \dots, \m_n$, which we will index
henceforth by $\m'_k$. Working over the quasi-local ring
\begin{equation}\label{Enakayama}
R_0 / I = R_0 / \cap_k \m'_k \simeq \times_k (R_0 / \m'_k),
\end{equation}
the images $M_i / (IM \cap M_i)$ are proper submodules of $M / IM$ (by
Nakayama's lemma) which also provide a covering. This shows that
$\sigma(M) \geq \sigma(M/IM)$, and we are reduced to the case where the
ring $R_0$ is a product of the fields $\F^{(k)} := R_0 / \m'_k$. Now
write $M = \oplus_k M^{(k)}$ with $M^{(k)} := \F^{(k)} M$; if
$\dim_{\F^{(k)}} M^{(k)} \leq 1\ \forall k$, then $M = R_0/I$ is cyclic
and cannot be covered by proper submodules, a contradiction. Thus there
exists $k$ such that $\dim_{\F^{(k)}} M^{(k)} \geq 2$, whence the vector
space $M^{(k)}$ admits a covering by proper $\F^{(k)}$-subspaces. The
proof is now completed by applying Proposition~\ref{Pproduct}. \qedhere
\end{enumerate}
\end{proof}

\subsection{Two settings involving covering finitely generated modules}

We now come to the case of possibly infinite minimal/irredundant
coverings of $R$-modules. Note that this case does not arise when
considering abelian groups, since in that case $R = \Z$ has no maximal
ideal $\m$ with infinite residue field $\m$. Thus, the results and
methods in this section necessarily lie outside the group covering
setting.

\begin{definition}
Let $R$ be a unital commutative ring, and $M$ an $R$-module admitting a
covering by proper $R$-submodules. Define the \textbf{covering number}
$\sigma(M)$ to be the smallest cardinal number $\aleph$ such that $M$ is
covered by $\aleph$-many proper submodules, but no fewer -- where `fewer'
means any cardinal that injects into $\aleph$ but is not in bijection
with it. (We use here the trichotomy of cardinal numbers, which is a
consequence of the Axiom of Choice; also, $\sigma(M)$ exists
e.g.~by~\cite{MeSa}.)
\end{definition}

In this subsection, we study coverings of finitely generated modules in
two settings. The first involves $R$ being a quasi-local ring -- i.e., a
commutative unital ring with only finitely many maximal ideals. In this
case, the same formula as in Theorem~\ref{Tfinite}(2) applies:

\begin{theorem}\label{Tquasilocal}
Let $R$ denote a quasi-local ring, and $M$ be a finitely
generated $R$-module. Then
\[
\sigma(M) = \nu'(R,M) = \min_{\m \in S} |R/\m|+1
\]
as cardinal numbers, where $\nu',S$ were defined in~\eqref{Eupperbound}.
(If $S$ is empty, then $M$ is cyclic.)
\end{theorem}

Thus, akin to Equation~\eqref{Enu1} and also the remarks following
Theorem~\ref{Tfinite}, once again the cardinal number $\nu'(R,M)$ turns
out to be the covering number.

As the arguments were essentially written out above, we merely sketch
this proof.

\begin{proof}
That $\sigma(M) \leq \min_{\m \in S} |R/\m|+1$ (when at least one such
$\m$ exists) is as in the proof of Theorem~\ref{Tfinite}(2). To show the
reverse inequality, let $M = \cup_{i \in \mathscr{I}} M_i$ be a covering
by proper submodules, where $|\mathscr{I}|$ is in bijection with
$\sigma(M)$. If $\m'_k$ denote the finitely many maximal ideals in $R$,
let $I = \cap_k \m'_k$ and work modulo $I$. Now the arguments
following~\eqref{Enakayama} apply verbatim to complete the proof.
\end{proof}

In our second setting here, the ring $R$ no longer needs to be
quasi-local, but contains a field $\F$ (with unit $1_R$).

\begin{theorem}
Suppose $R \supset \F \ni 1_R = 1_\F$ as above, and $M$ is a finitely
generated non-cyclic $R$-module.
Then $\sigma(M) \geq |\F| + 1$. If moreover $\F$ is an infinite field and
$R$ is countably generated over $\F$, then $\sigma(M) = |\F|+1$.
\end{theorem}

The point of interest here is that the (lower and upper) bounds are
universal for all finitely generated $R$-modules.

\begin{proof}
Take the smallest integer $n \geq 2$ such that $M$ has $n$ generators.
Now $R^n \twoheadrightarrow M$, so any covering of $M$ lifts to one of
$R^n$, whence $\sigma(M) \geq \sigma(R^n)$. Thus we reduce to $M = R^n$
with $n \geq 2$. Now suppose $R^n = \cup_{i \in \mathscr{I}} M_i$ is an
irredundant covering, so that each $M_i \neq R^n$ and $|\mathscr{I}| =
\sigma(R^n)$. Since $R \cdot \F^n = R^n$, it follows that $M_i \cap \F^n$
is a proper vector subspace of $\F^n$. But then,
\[
|\F|+1 = \sigma(\F^n) \leq |\mathscr{I}| = \sigma(R^n) \leq \sigma(M),
\]
where the first equality is by Proposition~\ref{Pvecunion}. This shows
the first assertion; the second will follow from the claim that
$\sigma(M) \leq |\F|+1$. To show the claim, first note that $|R| = |\F|$,
since every finitely generated $\F$-algebra has size $|\F|$. Now since
$M$ is finitely generated, $|M| = |\F|$; and since $M$ is not cyclic, it
is the union of its cyclic submodules, whence $\sigma(M) \leq |\F|$.
\end{proof}

In the setting of the preceding result, notice that every maximal ideal
$\m$ is an $\F$-vector subspace of $R$, whence every residue field $R/\m$
has the same size as $\F$. Given this and the earlier results in this
section, we end with the following natural question:\medskip

\textit{In what generality for the ring $R$, is it true that for all
finitely generated non-cyclic $R$-modules $M$, we have $\sigma(M) =
\nu'(R,M) = \min_{\m \in S} |R/\m|+1$ as in~\eqref{Eupperbound}?}

\section{Divisible groups and modules; direct sums of cyclic monoids}

Having understood covering numbers of abelian groups with finite
coverings, in particular finitely generated abelian groups, it is natural
to turn to other groups. In particular, every divisible abelian group
(such as $\mathbb{Q}$) is not a finite union but is a countable union of
proper subgroups. (The more general result -- over PIDs -- is mentioned
below.) The same result holds for groups of the form $\oplus_p (\Z / p
\Z)$, whenever the sum runs over an infinite set of pairwise distinct
prime integers. Indeed, that it is a countable union of proper subgroups
is as in the next paragraph, and it is not a finite union by
e.g.~Theorem~\ref{Traoreid}. This is an example of a $\Z$-module which is
not cyclic, but whose quotient modulo every maximal ideal is.

In a similar vein, for a field $\F$ it is easy to see (the assertion in
Proposition~\ref{Pvecunion}) that an infinite-dimensional vector space
over an infinite field is a countable but not finite union of proper
subspaces: simply write the basis as the union of a countable nested
sequence of nonempty proper subsets, and take their spans.

\begin{remark}\label{Rdvr}
There are also negative examples -- see e.g.~\cite[Example~2.7]{Zemlicka}
for the example of a local (in fact valuation) ring $(R,\m)$, with $M =
\m$ not a countable union of proper submodules.
\end{remark}

Yet another setting (which subsumes the case of abelian groups but not
vector spaces) is that of modules over a PID $R$. The case of finitely
generated torsion $R$-modules is easy to discuss (see later in this
section); for now, we mention the case of divisible modules.

\begin{prop}\label{Pdiv}
If $M$ is a nonzero divisible module over a PID $R$ (which is assumed to
not be a field), then $\sigma(M) = |\Z|$. More generally, $\sigma(M) \leq
|\Z|$ whenever $M$ is not reduced.
\end{prop}

\begin{proof}
No divisible module $M$ has a finite (in fact nonzero) quotient $M/\m M$,
hence is not a finite union of proper submodules by
Theorem~\ref{Tfinite}. On the other hand, by standard results -- see
e.g.~\cite[Exercises, \S 4.7]{Coh} and \cite{FuSa} -- every $R$-module is
the direct sum of a divisible (equivalently, injective) module and a
reduced module; and every divisible $R$-module is the direct sum of
copies of its field of fractions $\F$ and the `Pr\"ufer $p$-modules'
\[
M_p = R[p^\infty] := R[1/p]/R \subset \F/R, \qquad p \text{ prime } \in
{\rm Specm}(R).
\]
Now it is easy to show that the only $R$-submodules of $M_p$ are
$R[1/p^n]/R$, whence the assertion follows for $M = M_p$. Similarly, for
$M = \F$ (the quotient field of $R$), $\F$ is the nested union of
submodules $M'_n$, where we fix a prime $0 \neq p \in R$ and let
\[
M'_n := \{ a/b : a,b \in R, \ (a,b) = 1, \ p^n \nmid b \}, \qquad n \geq
1.
\]
Since at least one module from among $\F, M_p$ occurs in $M$, we have
$\sigma(M) = |\Z|$ for $M$ divisible. If instead $M$ is merely
non-reduced, then $M$ contains $\F$ or $M_p$ as a direct summand, whence
$\sigma(M) \leq |\Z|$.
\end{proof}

We conclude this part by discussing how far the above results take us, in
a restricted setting -- which in the above special cases (abelian groups,
vector spaces, or modules over PIDs) is a prominent case already studied
above. Namely, instead of working with finitely generated $R$-modules, we
instead consider direct sums of cyclic modules. In this case, the
`stopping point' seems to be the free module $M = R^2$ for general rings
$R$.

We now make this more precise, by writing down a sequence of
observations, and follow this by deducing several corollaries for
coverings of modules. We also apply these observations to minimally cover
direct sums of cyclic monoids.
Begin with a nonempty direct sum of \textit{nonzero} cyclic modules
\[
M = \oplus_{k \in \mathscr{K}} M_k
\]
over an arbitrary commutative unital ring $R$. 
\begin{enumerate}
\item Define $S$ and $\nu'(R,M)$ as in Theorem~\ref{Tfinite}(1). Then $S$
equals the collection of maximal ideals $\m$ that contain the
annihilators of $M_k$ for at least two $k \in \mathscr{K}$.

\item If $S$ is empty and $\mathscr{K}$ is finite, then $M$ is cyclic.
Indeed, if $M_k \cong R / I_k$ for all $k$, then each $I_k$ is contained
in exactly one maximal ideal, so the $I_k$ are coprime, and $M \cong R /
\cap_{k \in \mathscr{K}} I_k$ by the Chinese Remainder Theorem.

\item If $S$ is empty and $\mathscr{K}$ is infinite, then $M$ is covered
by a countable sequence of nested submodules (as discussed prior to
Remark~\ref{Rdvr}). Since $M$ is not a finite union of proper submodules
by Theorem~\ref{Tfinite}, we have $\sigma(M) = |\Z|$.

\item Otherwise, $S$ is nonempty. If $\nu'(R,M) < \infty$, then
Theorem~\ref{Tfinite} applies, and $\sigma(M) = \nu'(R,M)$.

\item Otherwise, $S$ is nonempty and $\nu'(R,M)$ is an infinite cardinal.
If $|\mathscr{K}| = \infty$, then $\sigma(M)$ is not finite (by
Theorem~\ref{Tfinite}), so $\sigma(M) = |\Z|$ (as discussed in a previous
case).

\item Otherwise, $S$ is nonempty, $|\mathscr{K}| < \infty$, and
$\nu'(R,M)$ is infinite.
If $R$ has finitely many maximal ideals -- i.e., $|{\rm Specm}(R)| <
\infty$, then $\sigma(M) = \nu'(R,M)$ by Theorem~\ref{Tquasilocal}.

\item (This is the `outstanding' case.) Otherwise, $S$ is nonempty,
$|\mathscr{K}| < \infty$, and $\nu'(R,M)$ and $|{\rm Specm}(R)|$ are both
infinite cardinals. As mentioned above, the torsion-free component of $M$
is the `stopping point' when one tries to understand minimal coverings.
\end{enumerate}

As a consequence of the above discussion, several special cases of rings
(and sums of cyclic modules over them) fall out as immediate corollaries.
We present here a sampling, starting with the final point above, where
the `torsion-free component' of $M$ is trivial, and we work over a
Dedekind domain:

\begin{cor}
Suppose $M$ is a finitely generated torsion module over a Dedekind
domain. Then $\sigma(M) = \nu'(R,M)$ if $S$ is nonempty, where
$\nu'(R,M)$ and $S$ are as in Theorem~\ref{Tfinite}(1); and if $S$ is
empty then $M$ is cyclic.
\end{cor}

\begin{proof}
It is well-known here that $M$ is a finite direct sum of cyclic torsion
$R$-modules, say $M_1, \dots, M_n$. Thus, we reduce the situation to
covering the finitely generated module $M / IM$ over the quasi-local ring
$R/I$, with $I = \cap_{k=1}^n \Ann_R (M_k)$. Now apply
Theorem~\ref{Tquasilocal}.
\end{proof}

The next corollary is immediate by Theorem~\ref{Tfinite}, and resolves
the above `outstanding' case under alternate additional assumptions on
$R$.

\begin{cor}
Suppose the assumptions in the final, outstanding case in the preceding
discussion apply: $M = \oplus_{k \in \mathscr{K}} M_k$ is a direct sum of
nonzero cyclic $R$-modules, where $S$ is nonempty, $|\mathscr{K}| <
\infty$, and $\nu'(R,M)$ and $|{\rm Specm}(R)|$ are both infinite
cardinals. If moreover $R$ is such that every residue field of $R/\m$ is
at most countable, then $\sigma(M) = |\Z|$, the smallest infinite
cardinal.
\end{cor}

As a third sample, we extend Proposition~\ref{Pvecunion} to the case of
local rings:

\begin{cor}\label{Clocal}
Suppose $(R,\m)$ is a local ring, and $M = \oplus_{k \in \mathscr{K}}
M_k$ is a direct sum of cyclic nonzero $R$-modules, with $|\mathscr{K}|
\geq 2$. Then $\sigma(M) = |\Z|$ if $R/\m$ and $\mathscr{K}$ are both
infinite, else $\sigma(M) = \nu'(R,M) = |R/\m|+1$.
\end{cor}

This too is shown by following the steps in the above discussion. There
is a similar result (and proof) for $R$ a direct product of finitely many
local rings.

\subsection{Covering sums of cyclic monoids}

We end with a final result, which is an application of the discussion
immediately above. We use the above results to answer a related variant:
that of covering monoids -- specifically, direct sums of cyclic monoids
-- in the spirit of finitely generated abelian groups or vector spaces.
Namely:\smallskip

\textit{Given a direct sum $M$ of cyclic monoids, how many proper
sub-monoids are required to cover $M$?}

\begin{theorem}[Sums of cyclic monoids]
Suppose $M$ is a direct sum of cyclic monoids. Then exactly one of the
following holds:
\begin{enumerate}
\item $M$ is a cyclic monoid, in which case it has no covering by proper
sub-monoids.
\item $M$ is an abelian group but not a cyclic monoid. In this case,
either the set $S$ defined in Theorem~\ref{Tfinite}(1) is empty and
$\sigma(M) = |\Z|$; or $S$ is nonempty, in which case $3 \leq \sigma(M) <
\infty$ and $\sigma(M) = \nu'(\Z,M)$ is obtained from
Theorem~\ref{Tfinite} (or Theorem~\ref{Traoreid}).
\item In all other cases, $\sigma(M) = 2$.
\end{enumerate}
\end{theorem}

In contrast, recall that no group is a union of two proper subgroups.

\begin{proof}
Suppose $M$ is not a cyclic monoid. If $M$ is a group and a direct sum of
cyclic monoids, then each factor is of the form $\Z / n_k \Z$ with $n_k
\geq 2$, for $k \in \mathscr{K}$, say. Now the above discussion,
following the proof of Proposition~\ref{Pdiv}, shows the assertion~(2) in
this corollary. Finally, if $M$ is not a group and not a cyclic monoid,
then one can write $M = M_1 \oplus M_2$, with $M_1 = \tangle{f_1}$ a
cyclic monoid that is not a group, and $M_2$ a nontrivial monoid. Then
$M$ is the union of $M_2$ and the monoid $(M \setminus M_2) \sqcup \{ 0
\}$; in fact this is a partition (in that the two sub-monoids intersect
only at $0$).
\end{proof}

\subsection*{Acknowledgments}
A.K.~was partially supported by Ramanujan Fellowship grant
SB/S2/RJN-121/2017, MATRICS grant MTR/2017/000295, and SwarnaJayanti
Fellowship grants SB/SJF/2019-20/14 and DST/SJF/MS/2019/3 from SERB and
DST (Govt.~of India), by grant F.510/25/CAS-II/2018(SAP-I) from UGC
(Govt.~of India), and by a Young Investigator Award from the Infosys
Foundation.



\end{document}